# A simulation-free seismic survey design by maximizing the spectral gap

Yijun Zhang[1], Mathias Louboutin[1], Ali Siahkoohi[1], Ziyi Yin[1], Rajiv Kumar[2], and Felix J. Herrmann[1]
[1] Georgia Institute of Technology
[2] Schlumberger

## Abstract

Due to the tremendous cost of seismic data acquisition, methods have been developed to reduce the amount of data acquired by designing optimal missing trace reconstruction algorithms. These technologies are designed to record as little data as possible in the field, while providing accurate wavefield reconstruction in the areas of the survey that are not recorded. This is achieved by designing randomized subsampling masks that allow for accurate wavefield reconstruction via matrix completion methods. Motivated by these recent results, we propose a simulation-free seismic survey design that aims at improving the quality of a given randomized subsampling using a simulated annealing algorithm that iteratively increases the spectral gap of the subsampling mask, a property recently linked to the quality of the reconstruction. We demonstrate that our proposed method improves the data reconstruction quality for a fixed subsampling rate on a realistic synthetic dataset.

## Introduction

Due to relatively recent breakthroughs in compressive sensing [Candès et al., 2006], seismic data is increasingly gathered randomly along spatial coordinates to decrease acquisition productivity [Mosher et al., 2014, Kumar et al., 2015, Chiu, 2019]. Wavefield reconstruction is used to recover fully sampled data from randomly subsampled observed seismic data [Hennenfent and Herrmann, 2008, Kumar et al., 2015, Zhang et al., 2020]. To remove the imprint of large gaps in uniform random sampling, Gilles [2008] proposed jittered subsampling, which by controlling the maximum gap size of subsampled data creates favorable conditions for seismic wavefield recovery based on sparsity promotion in a transformed domain made of localized atoms including curvelets [Herrmann et al., 2008]. While uniform random [Candès and Recht, 2009, Candès and Tao, 2010] and random jittered subsampling schemes [Herrmann and Hennenfent, 2008, Hennenfent and Herrmann, 2008] are relatively straightforward to generate, these sampling strategies are almost certainly suboptimal and have shown to be improvable by solving certain optimization problem [Mosher et al., 2014, Manohar et al., 2018, Li et al., 2017]. For instance, Mosher et al. [2014], Li et al. [2017], and later Titova et al. [2019] improved the reconstruction quality by devising a global optimization scheme that uses the mutual coherence.

In addition to wavefield reconstruction with optimized sampling schemes, Mosher et al. [2014] also proposed a simulation-based acquisition design to support the use of compressive sensing in seismic data acquisition. For time-lapse seismic, Guo and Sacchi [2020] also used a data-driven approach where the acquisition is optimized by using prior information on the seismic data [Manohar et al., 2018]. While these methods have lead to promising results, they either require significant computational resources to determine the optimal source-receiver layout using combined wavefield simulations and recoveries or require detailed information on the to-be-collected seismic data. Neither is feasible for the design of optimized sampling strategies in 3D.

To overcome this difficulty, this abstract provides a global optimization strategy for determining improved source-receiver layouts suitable for wavefield reconstructions based on matrix completion [Recht et al., 2010, Kumar et al., 2015, Kumar, 2017] without the need to carry out expensive wavefield simulations. Similar to Li et al. [2017] and Mosher et al. [2014], who propose a simulation-free optimization method based on the mathematical property of mutual coherence for transform-based wavefield reconstruction, our method involves improving the connectivity of graphs spanned by the binary sampling masks in the midpoint-offset domain. According to Bhojanapalli and Jain [2014], by maximizing the spectral gap—i.e., the distance between the two first singular values—of the binary sampling mask the connectivity of the graph is improved, which favors reconstruction by matrix completion, an observation recently confirmed by López et al. [2022] for 2D and 3D seismic wavefield reconstruction. While recent work by López et al. [2022] indeed negates the need to run multiple costly wavefield reconstructions for different candidate sampling masks, this work does not yet provide a constructive method to generate sampling masks that maximize the spectral gap.

Unfortunately, the design of acquisition masks that maximize the spectral gap is an NP-hard problem [Guo and Sacchi, 2020] whose solution requires a brute-force search through all combinatorial possibilities [Manohar et al., 2018]. When the number of sources or receivers becomes "large" [Li et al., 2016], this precludes its practical use; for example, there are $\binom{n}{m} = \frac{n!}{m!(n-m)!} = 75287520$ possible combinations when selecting $m = 5$ subsampling positions from a pool of $n = 100$ candidate sites. For this reason, we propose to obtain an approximate solution by maximizing the spectral gap using simulated annealing [Kirkpatrick et al., 1983], a stochastic local search optimization technique that is straightforward to implement, apply, and computationally feasible. The proposed method depends only on binary mask optimization, has a minimal computational cost, and should be adaptable to large-scale survey design.

We organize this expanded abstract as follows. First, we present the proposed optimization problem to maximize the spectral gap of subsampling masks. Next, we explain how to obtain the approximate acquisition masks via simulated annealing. We conclude by demonstrating numerical experiments on the 2D synthetic Compass dataset [Jones et al., 2012] and show an improvement in recovery quality compared to wavefield reconstruction of data collected with jittered subsampling method[Hennenfent and Herrmann, 2008].

## Methodology

Successful matrix-completion based seismic wavefield reconstruction [Recht et al., 2010, Kumar et al. [2015], Kumar [2017]] hinges on three critical factors, namely: (1) an appropriate randomized subsampling scheme, such as uniform random or jittered subsampling [Herrmann and Hennenfent, 2008, Hennenfent and Herrmann [2008]]; (2) existence of a transform domain in which the fully sampled seismic data organized as a matrix exhibits low-rank structure; (3) a computationally scalable matrix completion technique, which exploits the property that missing source and/or receivers increases the rank of these matrices. In this study, we propose a new constructive method to automatically generate improved source/receiver sampling masks, which favor seismic wavefield reconstruction via matrix completion in the midpoint-offset domain. We begin by describing our approach to acquisition design.

### Spectral ratio based acquisition design

Following López et al. [2022], we constitute the spectral gap by the spectral ratio (SR, the ratio of the first to second singular values), which becomes small for a large spectral gap. While the SR indeed has been shown to be a valuable quantity to predict the quality of wavefield reconstruction with matrix completion [López et al., 2022], ways to automatically generate acquisition masks with small SRs have so far been lacking. To meet this challenge, we cast the problem of finding optimized acquisition masks with small SRs as a minimization problem. Given $n_s$ source locations, $n_r$ receiver locations, and the source subsampling ratio $r$, we propose to solve a non-convex combinatorial optimization problem with respect to the subsampling mask

$\mathbf{M} \in \{0,1\}^{n_s \times n_r}$—i.e., we have

$$\mathcal{L}(\mathbf{M}) = \underset{\mathbf{M}}{\text{minimize}} \quad \frac{\sigma_2(\mathcal{S}(\mathbf{M}))}{\sigma_1(\mathcal{S}(\mathbf{M}))}$$
$$\text{subject to}$$
$$\|\mathbf{M}\|_0 = \lfloor n_s \times r \rfloor \times n_r \cap \mathbf{M} \in \mathcal{J} \cap \mathbf{M} \in \{0,1\}^{n_s \times n_r}. \quad (1)$$

In this optimization problem, the objective function consists of the spectral ratio (SR), defined by the ratio of the first, $\sigma_1(\cdot)$, and second, $\sigma_2(\cdot)$, singular values. $\mathcal{S}$ stands for the transformation operator with seismic reciprocity [Fenati and Rocca, 1984] from the source-receiver domain to the midpoint-offset domain. We constrain the solution to conserve the subsampling ratio ($\|\mathbf{M}\|_0 = \lfloor n_s \times r \rfloor \times n_r$) and to stay jittered sampled with $\mathcal{J} \subseteq \{0,1\}^{n_s \times n_r}$ being the set of all possible jittered subsampling acquisitions. This constraint guarantees that the spread of the survey will not be modified, but only the local source position will be optimized. The symbol $\lfloor \cdot \rfloor$ denotes the floor operation. As we previously mentioned, in order to solve this combinatorial optimization problem, we implemented a simulated annealing method to obtain a solution in a finite and acceptable time. We now describe this algorithm and link each step to its subsampling mask counterpart.

## Simulated annealing

Stochastic local search optimization algorithms are viable approximate methods for solving combinatorial optimization problems (e.g., Equation 1). Simulated annealing is a global optimization technique that uses local search to find approximate solutions to combinatorial optimization problems given a computational budget [Kirkpatrick et al., 1983, Şahin et al. [2010]].

This optimization method has three main components [Van Laarhoven and Aarts, 1987]: (1) an initial state, $\mathbf{M}_0$, representing the initial solution to the optimization problem (Equation 1); (2) a set of neighboring states for any given state, which will be used to update the current state randomly; and (3) a transition probability that determines the probability of moving from one state to another. During optimization, at each given state, $\mathbf{M}_k$, which represents the current solution to the optimization problem 1, a candidate state, $\tilde{\mathbf{M}}_k$, is chosen randomly from the neighboring states. Next, the algorithm transitions from the current state to the candidate state, i.e., from $\mathbf{M}_k$ to $\tilde{\mathbf{M}}_k$, if this transition reduces the objective function, i.e., $\mathcal{L}(\tilde{\mathbf{M}}_k) < \mathcal{L}(\mathbf{M}_k)$. On the other hand, if the objective function evaluated at the candidate state is larger than the current value, the algorithm makes the transition to the candidate state according to a transition probability, defined as follows [Kirkpatrick et al., 1983]:

$$p(\delta \mathcal{L}, k) = \exp\left(\frac{-\delta \mathcal{L}}{T(k)}\right), \quad (2)$$

where $k$ is the iteration number, $\delta \mathcal{L} = \mathcal{L}(\tilde{\mathbf{M}}_k) - \mathcal{L}(\mathbf{M}_k)$ indicates the change in the objective function (Equation 1) by moving to the candidate state, and $T(k) : \mathbb{R} \to \mathbb{R}^+$, typically called temperature function, is a monotonically decreasing function that reduces the uphill transition probability towards the end of optimization while allowing uphill movement early in the optimization. We choose the temperature function as $T(k) = T_0 \times \alpha^k$, following a geometric reduction rule, which is the most commonly used function in the simulated annealing with a start temperature $T_0$ and the decrease rate $\alpha$ [Kirkpatrick et al., 1983, Abramson et al. [1999]]. This allows the algorithm to escape from local minima in the initial stages of the optimization while ensuring downhill movement towards the end. Finally, the transition probability is smaller for candidate states that increase the objective function more, i.e., $\delta \mathcal{L} \gg 0$, minimizing the probability of moving to very bad solutions.

To adapt simulated annealing to the acquisition design optimization problem (cf. Equation 1), we define the states as arbitrary positioning of sources. The algorithm is summarized in Algorithm 1. This algorithm is initialized with a subsampling mask $\mathbf{M}_0$ that is generated by using jittered subsampling method, known to facilitate seismic wavefield recovery [Hennenfent and Herrmann, 2008]. After updating the temperature function $T(k)$ (line 1) [Ma, 2002], we select a source position $\tilde{\mathbf{M}}_k$ within the neighborhood of the current position (line 2). We then update the source position to this updated state according to the loss decrease and probabilistic update rule (lines $3-8$). After a predetermined number of iterations, the algorithm outputs the source sampling mask $\mathbf{M}_K$ with smaller SR.

---

**Algorithm 1** Spectral ratio minimization via simulated annealing.

**Inputs:**
    $\mathbf{M}_0$: Initial source positions using jittered subsampling method.
    $K$: Maximum number of iterations.
    $T(k)$: Temperature function.
    $p$: Transition probability (cf. Equation 2).

0. **for** $k = 0$ **to** $K - 1$ **do**
1.     $T(k) = T_0 \times \alpha^k$
2.     $\tilde{\mathbf{M}}_k \leftarrow$ randomly pick a neighboring state
3.     $\delta \mathcal{L} = \mathcal{L}(\tilde{\mathbf{M}}_k) - \mathcal{L}(\mathbf{M}_k)$
4.     **if** $\delta \mathcal{L} \leq 0$
5.         $\mathbf{M}_{k+1} = \tilde{\mathbf{M}}_k$
6.     **else**
7.         $\mathbf{M}_{k+1} = \begin{cases} \tilde{\mathbf{M}}_k & \text{with probability } p(\delta \mathcal{L}, k) \\ \mathbf{M}_k \end{cases}$
8.     **end if**
9. **end for**

**Output:** $\mathbf{M}_K$

---

In order to satisfy the constraints introduced in Equation 1, we carefully define the neighborhood of acceptable state to prevent sources positions to cluster around a few areas of the survey. We now detail its design.

**Neighboring states**

Given the source sampling at an iteration, we randomly select 20% of the source positions in the current state to balance between exploring the search space and avoiding too large change between adjacent iterations [Olorunda and Engelbrecht, 2008, Assad and Deep [2018]]. Using the subsampling factor $f = \frac{1}{r}$, the fine grid with all possible source locations is divided into $f$ equal regions. Each selected source is allowed to randomly shift within the region in which it is located. For clarity, we summarize the perturbation rule in Figure 1 with $n_s = 20$ source positions and a subsampling factor of $f = 5$. We choose the movement range $R$ to ensure that we remain close to the jittered sampling, which has been shown to result in better wavefield recoveries [Hennenfent and Herrmann, 2008]. We now detailed our synthetic numerical experiment demonstrating the benefits of our method for data reconstruction.

# Numerical experiments

We consider a 2D marine dataset simulated over the realistic Compass model [Jones et al., 2012]. The dataset consists of 300 sources and 150 receivers sampled at 12.5 m. The data is recorded at a 2 ms sampling rate for 2.046 s (1024 time samples). Based on this dataset, we proceed in two steps. First, we will compare the subsampling mask we obtain with our proposed method against the standard jittered sampling mask. Second, we will show that the recovered data is of better quality as expected from the optimal SR that quantifies the expected quality of recovery. The jittered subsampling method used in this abstract allows neighboring subsamples (non-gap subsamples), which creates favorable conditions for recovery and is defined as optimally-jittered subsampling in Hennenfent and Herrmann [2008]. We use the weighted matrix completion method [Zhang et al., 2020] to recover the observed data and evaluate the quality of the recovered dataset.

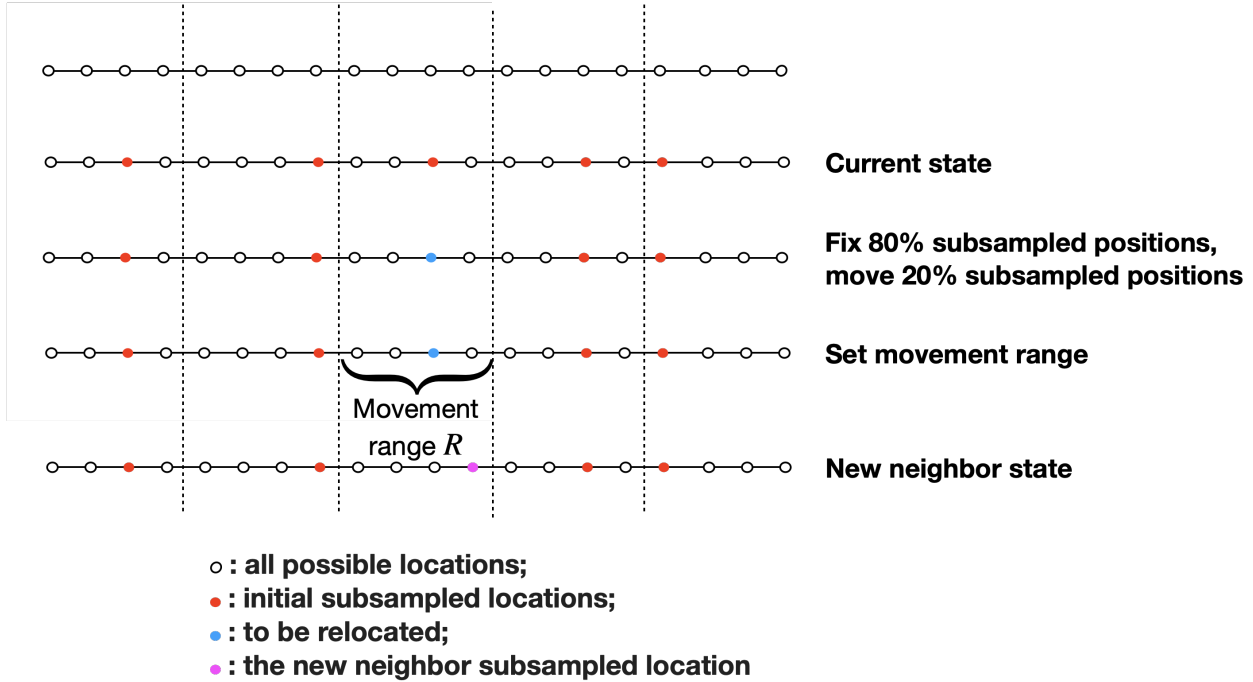

○ : all possible locations;
● : initial subsampled locations;
● : to be relocated;
● : the new neighbor subsampled location

Figure 1: Random state perturbation rules to satisfy the constraints. White circles (○) indicate all possible source locations. Five red circles represent an initial state. The blue circle represents the 20% of sources that we will move within movement range $R$ to arrive at a neighboring state (purple circle).

We start with picking a subsampling mask using the jittered method [Herrmann and Hennenfent, 2008] that includes 20% of the sources. In Figure 2a, we show the subsampling mask in the source-receiver domain. Under the assumption of the source-receiver reciprocity [Fenati and Rocca, 1984], we apply this reciprocity on the subsampling mask to implement a realistic seismic survey design. The subsampling mask in the midpoint-offset domain with seismic reciprocity is depicted in Figure 2b. With this mask as an initial guess, we perform 4000 iterations of simulated annealing to obtain an optimal subsampling mask. Figure 2c and Figure 2d show the resulting subsampling mask in the source-receiver and midpoint-offset domains, respectively. We observe that the SR was reduced by 30% with a fixed subsampling rate hinting towards well-improved data reconstruction. The proposed method is a simulation-free method that depends only on binary mask optimization.

With this optimized subsampling mask, we now perform data reconstruction via weighted matrix completion [Zhang et al., 2020] and compare the result against reconstructing the data sampled with the initial jittered subsampled mask. In both cases, we use the same algorithm and hyperparameters (e.g., number of iterations, rank) for a fair comparison. We summarize the recovery in Figure 3.

We first show the ground truth in Figure 3a, where the right plot shows the full shot record and the left one depicts the later arrival events between about 1 s to 2 s. By applying these two masks (jittered mask and proposed mask) individually on the ground truth, we obtain two observed datasets with 80% of sources missing. The proposed subsampled data is illustrated in Figure 3b. Figure 3c shows the reconstruction from jittered observed data with a signal-to-noise ratio (SNR) of 14.6 dB for the full shot record and 12.8 dB for the later arrival events. The reconstruction from the proposed subsampled data is shown in Figure 3d, with SNRs of 14.91 dB for the full shot record and 12.9 dB for the later arrival events. Figure 3e illustrates the difference between Figure 3c and Figure 3a, whereas Figure 3f shows the difference between Figure 3d and Figure 3a. The wavefield reconstructions demonstrate that the reconstruction from the proposed subsampling mask gives a more accurate data reconstruction for the full shot record and the later arrival events in terms of SNR. The difference is significantly reduced in Figure 3f in contrast to Figure 3e.

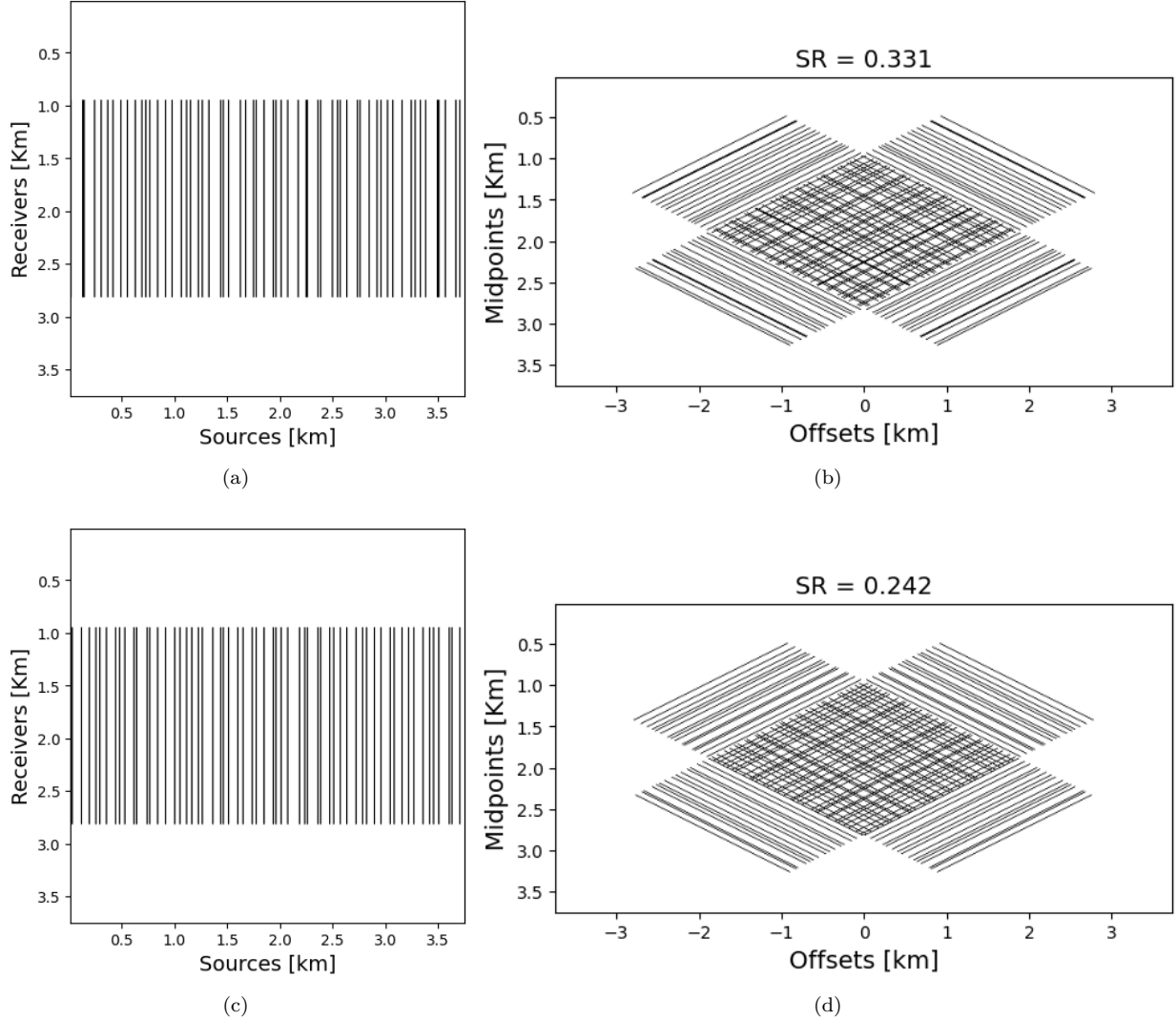

Figure 2: Jittered subsampling mask in the *(a)* source-receiver domain and *(b)* midpoint-offset domain (SR = 0.331). Optimized subsampling mask in the *(c)* source-receiver domain and *(d)* midpoint-offset domain (SR = 0.242).

To further validate the performance of our proposed method, we show that our proposed mask outperforms on the average standard jittered acquisition and not just for a single experiment. We randomly generate five independent jittered subsampling masks (removing 80% of sources) and then utilize the proposed approach to minimize the SR of these five jittered subsampling layouts. These five jittered and proposed masks are then used to perform weighted matrix completion [Zhang et al., 2019, Zhang et al. [2020]] and reconstruct the full dataset. The results are summarized in Figure 4. The bar plots in Figure 4 lead to the following observations. First, our proposed method consistently reduced the spectral ratio by at least 11% leading to a similar optimal SR for this given subsampling ratio and acquisition. Second, the recovered data always presents a higher SNR representative of a more accurate wavefield reconstruction. These two results show that despite being a potentially aleatory method, our simulated annealing based SR minimization method consistently provides a subsampling mask best fitted for data reconstruction.

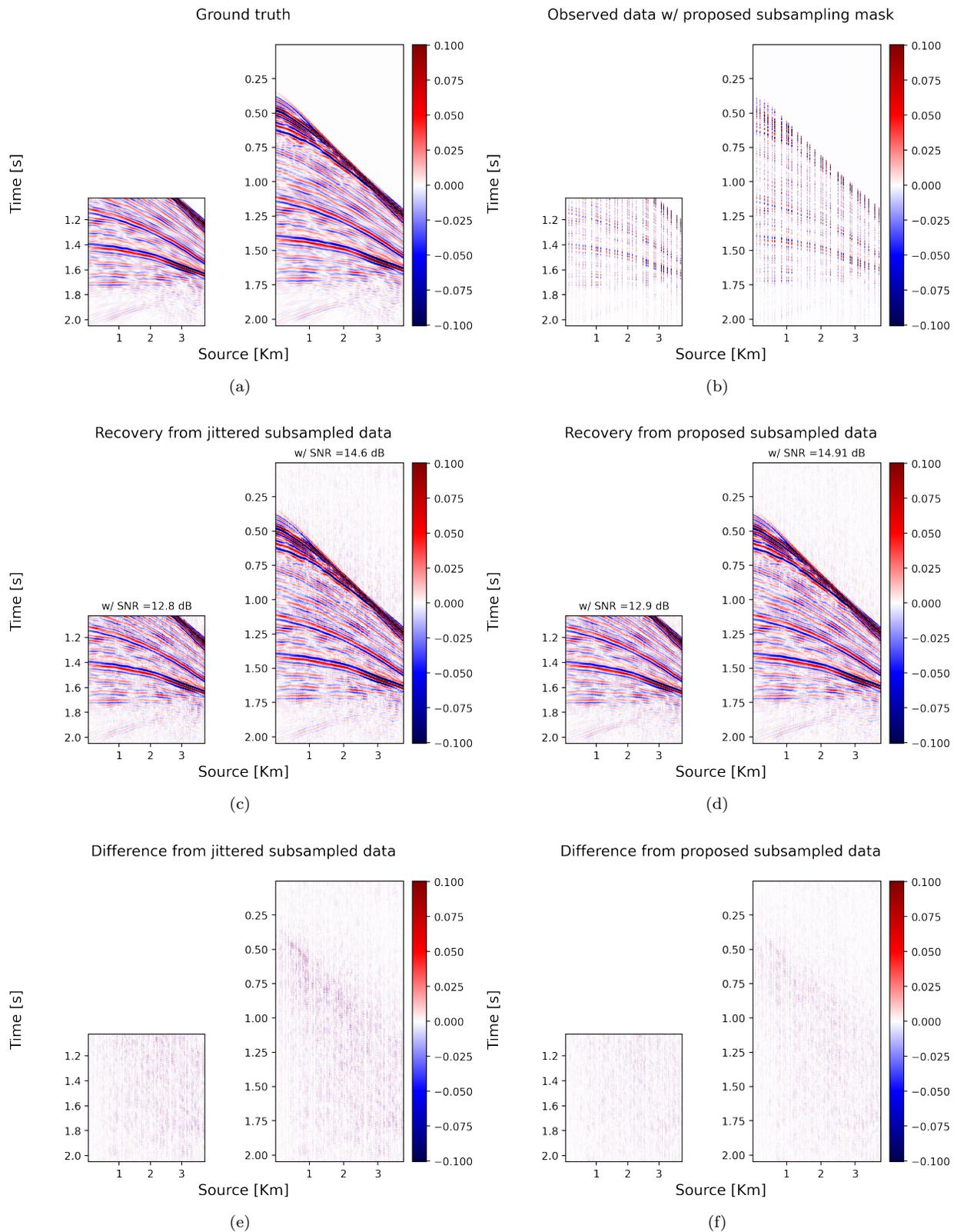

Figure 3: Wavefield reconstruction results in the time domain. *(a)* Ground truth. *(b)* 80% subsampled seismic data with proposed subsampling. Reconstructions from 80% missing sources: *(c)* jittered subsampling, SNR = 14.6 dB and 12.8 dB for later arrival events, *(d)* improved subsampling with SNR = 14.91 dB and 12.9 dB for later arrival events.

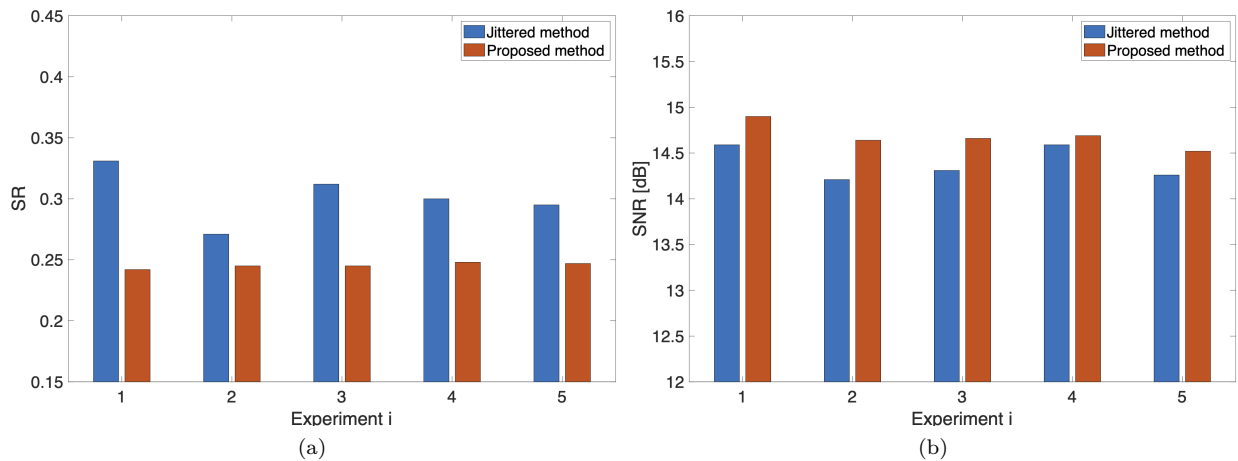

Figure 4: *(a)* SR comparison (lower is better) of subsampling masks using jittered method versus proposed method. *(b)* Reconstruction SNR comparison (higher is better) from observed data using jittered method and proposed method. The results are obtained by five independent experiments.

## Conclusions

We proposed a simulation-free method for seismic survey design in this study by minimizing the spectral ratio using the simulated annealing algorithm. Because the proposed method solely relies on a binary mask optimization rather than being data-driven, the computational cost is minimal and should scale to industry-size survey design. Through analysis and experiments, we conclude that the proposed method leads to an optimal subsampling mask best fitted for wavefield reconstruction based on matrix completion. Future work will focus on applying the proposed method to the three-dimensional field data.

## Acknowledgement

This research was carried out with the support of Georgia Research Alliance and partners of the ML4Seismic Center.